\documentclass[11pt]{article}
\catcode`\@=11 \catcode`\@=12

\usepackage{amsmath}
\usepackage{amssymb}
\usepackage{amsthm}
\usepackage{oldgerm}
\usepackage{multicol}

\renewcommand{\d}{{\rm{d}}}

\newcommand{\Rm}{\mathbb{R}}

\newcommand{\mL}{\mathcal{L}}

\newcommand{\mC}{\ensuremath{\mathcal{C}}}

\newcommand{\mG}{\ensuremath{\mathcal{G}}}

\newcommand{\mY}{\ensuremath{\mathcal{Y}}}
\newcommand{\mN}{\ensuremath{\mathcal{N}}}

\newcommand{\ca}{\ensuremath{\text{cart}_k}}

\newcommand{\mP}{\ensuremath{\mathcal{P}}}

\newcommand{\vs}{\vspace{.2cm}}
\newcommand{\J}{\ensuremath{J^1(N,M)}}

\newtheorem{lem}{Lemma}
\newtheorem{thm}{Theorem}

\newtheorem{cor}[lem]{Corollary}
\newtheorem{prop}[lem]{Proposition}
\newtheorem{defn}[lem]{Definition}

\def\proof {\noindent{\sc{Proof. }}}
\def\qed {\mbox{}\hfill {\small \fbox{}} \\}
\def\lto{\longrightarrow}
\def\lmto{\longmapsto}

\def\leq{\leqslant}
\def\geq{\geqslant}

\title{Young measures, Cartesian maps, and polyconvexity}
\author{Patrick  Bernard\footnote{membre de l'IUF} 
\hspace{.2cm}and Ugo Bessi}
\date{April 2008}
\begin{document}

\maketitle 

\vspace{1cm}
\begin{center}
-----
\end{center}
\begin{multicols}{2}

\begin{small}
\noindent
Patrick Bernard,
 Universit\'e  Paris-Dauphine,\\
CEREMADE, UMR CNRS 7534\\
Pl. du Mar\'echal de Lattre de Tassigny\\
75775 Paris Cedex 16,
France\\
\texttt{patrick.bernard@ceremade.dauphine.fr}\\

\noindent
Ugo Bessi, Universit\`a\ di Roma 3\\
Dipartimento di Matematica\\
Largo S. Leonardo Murialdo 1\\
00146 Roma, Italia\\
\texttt{bessi@matrm3.mat.uniroma3.it}\\

\end{small}

\end{multicols}
\vs
\thispagestyle{empty}
\begin{center}
-----
\end{center}

\textbf{Abstract. } 
We consider the variational problem consisting of minimizing a
polyconvex integrand for maps between manifolds. We offer a
simple and direct proof of the existence of a minimizing map.
The proof is based on Young measures.\\
\begin{center}
-----

\end{center}

\textbf{R\'esum\'e. }
On consid\`ere le probl\`eme variationnel consistant \`a minimiser
un int\'egrand polyconvexe pour des applications entre vari\'et\'es.
On donne une preuve simple et directe de l'existence d'une application minimisante en utilisant les mesures d'Young.\\
\begin{center}
-----
\end{center}

\textbf{Riassunto.} 
Consideriamo il problema di minimizzare un funzionale policonvesso tra variet\`a. Diamo una dimostrazione semplice e diretta dell'esistenza di un minimo usando le misure di Young.
\begin{center}
-----
\end{center}


\newpage
\section{Introduction}
Let $\bar N$ be a compact Riemannian manifold 
with boundary of dimension
$n$ 
and $M$ be a complete Riemannian manifold of dimension  $m$. 
We denote by $N$ the interior of $\bar N$,
and by $\d t$ the non-negative Borel  measure on $N$ associated
with the metric. It can be defined as the $n$-Hausdorff measure 
associated with the Riemannian distance, and it exists also if $N$ is not orientable.
This measure is also characterized by the fact that 
$$\int _D \d t=\left |\int_D \Omega\right|$$
for each embedded disk $D\subset M$, where $\Omega$ is any of the 
two unitary volume forms on $D$.
We set $p=\min\{m,n\}$.
We want to discuss the classical question of minimizing 
the quantity 
$$\int_N L(t,u(t),\d u_t) \;\d t
$$
on appropriate classes of maps $u:N\lto M$. 
We recall some terminology about Lagrangians:\vs

\noindent
\textbf{Convexity :}
Given $k\in \{1,\ldots,\min(m,n)\}$, we say that the integrand 
$L(t,x,v)$ 
is $k$-convex if it 
can be written in the form
$$
L(t,x,v)=\textsf{L}(t,x,v,\wedge_2 v,\ldots,\wedge_k v)
=\textsf{L}_{t,x}(v,\wedge_2 v,\ldots,\wedge_k v)
$$
with a Borel measurable function $\textsf{L}$ such that
$$\textsf{L}_{t,x}
:\mL(T_tN,T_xM)\times \mL(\wedge_2 T_tN,\wedge_2T_xM)\cdots
\times \mL(\wedge_kT_tN,\wedge_k T_xM))
\lto \Rm\cup \{+\infty\}
$$
is convex for each $t$ and $x$.
When $k=1$, this just says that $L$ is convex in $v$;
when $k=p$, this hypothesis is usually 
called polyconvexity.
Let us explain the notations.  We denote by $\mL(E,F)$
the set of linear maps between $E$ and $F$. The space $\wedge_iE$
is the vector space $E\wedge E\ldots \wedge E$
generated by elements of the form $e_1\wedge e_2\ldots \wedge e_i$,
$e_j\in E$.
Given $v\in \mL(E,F)$, we denote by 
$\wedge _i v\in \mL(\wedge_iE,\wedge _iF)$ the linear map such that 
$$
\wedge_iv(e_1\wedge \ldots\wedge e_i)=
v(e_1)\wedge v(e_2) \ldots \wedge v(e_i).
$$ 
If $E$ has dimension $n$ and $F$ has dimension $m$,
the dimension of $\wedge _i E$ is  $C^i_n$ (binomial coefficient)
and $\wedge_i v$ can be represented by a $C^i_n\times C^i_m$
matrix whose coefficients are the determinants of all $i\times i$
sub-matrices that can be extracted from the matrix representing $v$.
Let us denote by $J^1(N,M)$ the manifold of $1$-jets of maps from
$N$ to $M$. If  $N$ is an open subset of $\Rm^n$ and $M=\Rm^m$,
then
$$J^1(N,M)=N\times M\times \mL(\Rm^n,\Rm^m).$$ 

\noindent
\textbf{Regularity :}
We say that $L:\J\lto \Rm\cup \{+\infty\}$ is a normal integrand if 
it is a Borel measurable function and if, 
for almost all $t\in N$,
the function $(x,v)\lmto L(t,x,v)$ is lower semi-continuous.\vs

\noindent
\textbf{Coercivity :}
We say that $L$ is $k$-superlinear if there exists
a superlinear function $l:[0,\infty)\lto \Rm$
such that 
$$
L(t,x,v)\geq l(\|v\|+\|\wedge_2 v\|
+\cdots +\|\wedge_k v\|)
$$
for all $(t,x,v)\in J^1(N,M)$.\vs

Our goal  is to explain
a compact and simple approach to this kind of problems.
We do not present real novelties, and much of
the techniques we will use can be found in \cite{Pe}  or 
\cite{GMS1,GMS2}.
Yet we  believe it is not useless to present   the short path to Theorem \ref{main}
that follows. This work started with an attempt to extend the methods
of \cite{BY} to higher dimension.

We will define, by studying relevant sets of Young measures,
  sets $\ca(N,M)$ of maps $u:N\lto M$
such that 
$$W^{1,n}(N,M)\subset \text{cart}_{\min\{n,m\}}(N,M)\subset
 \ca(N,M)
 \subset
 \text{cart}_1(N,M)=  W^{1,1}(N,M).$$
In the case $k=\min\{n,m\}$ the set that we denote $\ca(N,M)$
is similar to the set denoted 
$\text{cart}^1(N,M)$  in \cite{GMS1,GMS2},
 but our presentation is quite different.

In order to state appropriately a variational
problem, it is useful to specify boundary conditions.
We assume that the boundary $\partial N$  of $\bar N$ is not empty
and we fix a map $u_0\in \ca(N,M)$.
We denote by $\ca(N,M;u_0)$ the set of maps $u$ 
in $\ca(N,M)$
such that the trace of $u$ in $L^1(\partial N, M)$
is equal to the trace of $u_0$. These traces are well defined
(at least in the case where $M$ is a Euclidean space)
because $\ca(N,M)\subset W^{1,1}(N,M)$, and each element of 
$W^{1,1}(N,M)$ has a unique boundary trace in $L^1$ (see for instance $\cite{EG}$ for the definition); we recall that the integration by parts formula holds for this trace.
In the case where $M$ is a manifold, we shall give the 
precise definition of $\ca(N,M;u_0)$ in Section \ref{compactsect}.
Our goal is to provide a short and direct proof of the following result:

\begin{thm}\label{main}
Let $k\in \{1,\ldots,\min(m,n)\}$ be given
and let $L:\J\lto \Rm\cup \{+\infty\}$  be a $k$-convex and $k$-superlinear 
normal integrand.
Let $u_0\in \ca( N,M)$ be given,
such that 
$$
\int _N L(t,u_0(t), \d u_0 (t)) \d t<\infty.
$$
Then there exists a map $u\in \ca(N,M;u_0)$
which minimizes the integral
$$
\int _N L(t,u(t), \d u(t))\d t
$$
in $\ca(N,M;u_0)$.
\end{thm}

Note that, in general,   we may 
have 
$$
\inf _{u\in W^{1,1}(N,M;u_0)} \int_N L(t,u(t),\d u_t)\d t
<
\min _{u\in \ca(N,M;u_0)}\int_N L(t,u(t), \d u_t)dt
$$
and, even if $u_0$ is smooth, 
$$
\min _{u\in \ca(N,M;u_0)}\int_N L(t,u(t),\d u_t) \d t
<
\inf _{u\in C^1(N,M;u_0)} \int_N L(t,u(t),\d u_t)\d t.
$$
The paper \cite{BM} contains an example in which $N$ has dimension 1 and the minimum on 
$W^{1,1}(N,\Rm^n)=\text{cart}_1(N,\Rm^n;u_0)$ 
is smaller than the infimum on $C^1$.
As for the other strict inequality, we are going to see an example in section \ref{example}.

Theorem \ref{main} is a slight extension on the seminal results
of Ball, \cite{Ba:77}.
Compared to this work and to classical papers on polyconvexity,
our proof works under slightly weaker coercivity. Our Theorem reduces to the famous Tonelli theorem
in the case $n=1$, and to the famous De Giorgi Theorem in the case
$m=1$.
 Several extensions are known, which go much beyond
 what we plan to expose. First, the kind
of convexity hypothesis can be relaxed to the so-called quasi-convexity,
but then one has to add more stringent growth conditions,
see \cite{Da,Pe,Bu}.
Second, one can, in certain circumstances, relax the
coercivity condition to the case when $L$ has only linear growth,
by using cartesian currents and functions of bounded variations.
Excellent surveys of these methods are in \cite{GMS1,GMS2,Bu}.

Our approach is based on Young measures, also called
parametrized measures. A survey on the use of Young measure
is the book \cite{Pe}.
Many of our techniques are adapted from this book.
 In section \ref{generalized}, we define the various
 sets of Young measures that are useful, in particular the
 set of Generalized Maps, on which it is appropriate
 to relax the variational problem.
We study the structure of generalized maps and conclude 
that minimizing generalized maps correspond to minimizing maps.
In section \ref{compactsect}, we prove a compactness results 
under boundary conditions.
In section \ref{conclusion}, we briefly expose how the various tools
exposed in sections \ref{generalized} and \ref{compactsect}
lead to a proof of Theorem \ref{main}.
We also collect various related remarks.
Finally, we recall in the Appendix some relevant facts on the 
topology of some spaces of measures.

We end this introductory section collecting some notation and material on $n$-vectors and $n$-forms we shall need in the following.

\subsection{Some algebra}

Let $E$ and $F$ be two  Euclidean vector spaces
of dimension $n$ and $m$. It will be convenient to set $p:=\min(m,n)$.
We denote by $\mL(E,F)$ the set of linear maps
between $E$ and $F$.
Recall that the vector space  $\wedge _l E$ is endowed with
a natural inner product (which is induced from the inner product 
of $E$).
This inner product can be characterized by the property that 
$$
\langle v_1\wedge\ldots \wedge v_l,w_1\wedge \ldots
\wedge w_l\rangle=
\det(G)
$$
where $G\in M^{l,l}(\Rm)$ is the Gram
matrix $G_{i,j}=\langle v_i,w_j\rangle$
and $\det G$ is the determinant of $G$.
Note  that 
$$\|v_1\wedge \ldots \wedge v_l\|
:=
\langle v_1\wedge \ldots \wedge v_l,v_1\wedge \ldots \wedge v_l
\rangle ^{1/2}=1$$
if $(v_1,\ldots v_l)$ is orthonormal in $E$.
Given $a\in \mL(E,F)$, we denote by 
$\|a\|:= \sup_{\|x\|_E\leq 1} \|a(x)\|_F$ its norm and by
$\wedge_l a\in\mL(\wedge _l E,\wedge_l F)$
the unique linear map which satisfies
$$
\wedge_k a (v_1\wedge v_2\ldots \wedge v_l)
=a(v_1)\wedge a(v_2)\wedge \ldots \wedge a(v_l)
$$
for each $v_1,\ldots, v_l$ in $E$.
The map $\wedge_l a$ is called the $l$-adjoint of $a$.
In coordinates, this map is represented by a matrix
whose coefficients are the $l$-minors of $a$.

If $\omega$ is a $k$-form on $E$ and $U$ is an $l$-vector,
$l\leq k$, we denote by $i_U\omega$
the $(k-l)$-form defined by 
$$i_U\omega\cdot v=\omega\cdot(U\wedge v)$$
for any $(k-l)$-vector $v$.

\begin{lem}\label{calcul}
Let $\Omega$ be a volume form on $E$, let $\lambda$ be the unique 
$n$-vector on $E$ such that $\Omega\cdot \lambda=1$,
let $U$ be a $k$-vector on $E$ and let $a\in \mL(E,F)$. 
Then, for any $k$-form $\chi$ on $F$, we have 
$$
\big((i_U\Omega)\wedge \chi\big)\circ \wedge_n(Id\oplus a)
\cdot \lambda
=
(-1)^{k(n-k)}\chi\circ \wedge_ka\cdot U,
$$
where $Id\oplus a:E\lto E\oplus  F$ is the map $v\lmto (v,a(v))$.
\end{lem}
\proof
We make a proof in coordinates.
Let $(e_1,\ldots e_n)$ be a base of $E$ such that 
$\lambda=e_1\wedge e_2\wedge \ldots\wedge e_n$.
If $e_j^*$ is the dual base of $E^*$, then we have
$\Omega=e_1^*\wedge\ldots\wedge e_n^*$.
If $I$ is a subset of $\{1,\ldots, n\}$,
we denote by $e_I$ the product 
$e_{\alpha_1}\wedge \ldots \wedge e_{\alpha_i}$,
where $i$ is the cardinal of $I$, and $\alpha_j, 1\leq j\leq i$
are the elements of $I$ in increasing order.
We denote by $\sigma(I)$ the sign such that 
$e_{I}\wedge e_{I^c}=\sigma(I)\lambda$, where 
$I^c$ is the complement of $I$.
Note that $\Omega=\sigma(I)e_I^*\wedge e_{I^c}^*$,
so that   $i_{e_I}\Omega=\sigma(I)e^*_{I^c}$.
Note that $\sigma(I^c)=(-1)^{k(n-k)}\sigma(I)$.
It is sufficient to prove the Lemma for 
$U=\sigma(J^c)e_{J^c}$, where $J$ has cardinal $n-k$,
in which case $i_U\Omega=e^*_J$.
We have
$$
\wedge_n(Id\oplus a)\cdot \lambda=\sum_I \sigma(I)e_I
\wedge (\wedge _{|I^c|}a\cdot e_{I^c})
$$
where the sum is taken on all subsets $I$ of 
$\{1,\ldots,n\}$ and $|I|$ is the cardinal of $I$.
We get 
$$
\big((i_U\Omega)\wedge \chi\big)\circ \wedge_n(Id\oplus a)
\cdot \lambda
=
(e^*_J\wedge \chi)\circ \wedge_n(Id\oplus a)\cdot\lambda
=
$$
$$\sigma(J)\chi\circ \wedge_ka\cdot e_{J^c}
=(-1)^{k(n-k)}\chi \circ \wedge_k a \cdot U.
$$

\qed

\section{Generalized maps}\label{generalized}

Let us first recall the definitions of the Sobolev space 
$W^{1,q}(N,M)$. We say that $u$ belongs to this space if there
exists a $v(t)\in \mL(T_tN,T_{u(t)}M)$, depending measurably on $t$,
such that 
$$\int_N\|v(t)\|^q_{t,u(t)} \d t<\infty
$$
and
$$\d (\chi\circ u)_t=\d\chi_{u(t)}\circ v(t)$$
in the sense of distributions for all bounded smooth functions
 $\chi:M\lto \Rm$ with bounded derivative.
This can be written intrinsically on the manifold $N$ by requiring
that, for each smooth vectorfield $U(t)$ on $N$ compactly supported
in the interior of $N$,
we have
$$
\int_N \d \chi_{u(t)}\circ v(t)\cdot U(t) \d t+
\int_N \chi(u(t))\cdot \text{\rm div}U(t)\d t=0.
$$
 It is not hard to see that, if $M$ is a Riemannian submanifold
of a Euclidean space $E$, then $W^{1,q}(N,M)$ is just the set of
 the functions $u\in W^{1,q}(N,E)$ which satisfy $u(t)\in M$
for almost every $t$.
We recall that smooth functions are not
necessarily  dense in these spaces
if $q<n$.

\subsection{Young measures}
Let us denote by $\J$ the set of
$1$-jets of maps $u:N\lto M$.
In many examples, $N$ is an open subset of $\Rm^n$,
$M$ is $\Rm^m$, and then 
$$
\J=N\times \Rm^m\times \mL(\Rm^n,\Rm^m).
$$
We shall usually denote by $(t,x,v)$ the points of $\J$.
We define the function 
$$
r_k(t,x,v)=1+\|v\|+\|\wedge_2 v\|+\cdots+\|\wedge_k v\|
$$
and associate to it the complete metric space $\mP_{r_k}(\J)$
as in the appendix.
This is the space of Borel probability measures 
$\eta$ on $\J$
such that $\int r_k d\eta<\infty$. 
We note that the measure on $N$ induced by the Riemann metric, which we have denoted by $\d t$, is finite, since $N$ is compact; to simplify the following definitions, we shall suppose that the measure of $N$ is $1$. 

Let $t:\J\lto N$ denote the natural projection; 
we denote by $\mY_k(N,M)$
the set of non-negative Borel measures 
$\eta\in \mP_{r_k}(\J)$ such that 
$t_{\sharp}\eta$ coincides with the  measure $\d t$.
We endow $\mY_k(N,M)$ with the topology induced from 
$\mP_{r_k}(\J)$.

\begin{prop}\label{lsc}
Assume that  $L$ is a normal integrand 
which is bounded from below
(or more generally such that $L/r_k$
is bounded from below),
then $\eta\lmto \int L d\eta$ is lower semi-continuous on 
$\mY_k$.
\end{prop}

\proof
Assume first that $L$ is continuous and that $L/r_k$ is bounded.
Then, the functional is continuous by definition of the topology on
$\mP_{r_k}$.

As an intermediate step, assume that 
$L(t,x,v)$ is a Caratheodory integrand
(measurable in $t$ and continuous in $(x,v)$) and that $L/r_k$
is bounded. 
By the Scorza-Dragoni Theorem, (see 
\cite{BeLa:73}, Theorem I.1.1, p 132)
there exists an increasing sequence $K_i$ of compact subsets on $N$
such that $L$ is continuous on $J^1(N,M)_{|K_i}$
(the set of points $(t,x,v)$ such that $t\in K_i$)  and such that 
$\cup_i (K_i)$ has full measure in $N$.
Then, there exists a sequence of continuous functions
$L_i$ such that $|L_i|/r_k$ is bounded,
independently of $i$, and such that $L_i=L$
on $J^1(N,M)_{|K_i}$.
It follows that the map $\eta\lmto \int L\d \eta$ 
is the uniform limit on $\mY_k(N,M)$
of the continuous maps $\eta\lmto \int L_i\d \eta$,
and therefore it is continuous on $\mY_k(N,M)$.

In the general case,
we first write the integrand $L(t,x,v)=r_k(t,x,v)g(t,x,v)$
with a normal integrand $g$ which is bounded from below. 
Then $g$ is  the increasing pointwise limit
of a sequence $g_i$  of bounded Caratheodory integrands,
see \cite{BeLa:73}, Theorem I.1.2, p 138.
Finally,  the map $\eta\lmto \int Ld\eta$ is the increasing limit
of the continuous maps $\eta\lmto \int r_kg_i  \d\eta$, and therefore it 
is lower semi-continuous.

\qed

\subsection{Closed measures}
It is a fundamental and well known observation that
there exists many null-Lagrangians, that is functions
$F(t,x,v):J^1(N,M)\lto \Rm$ such that 
$$\int_N F(t,u(t),\d u_t)\d t=0
$$
for all $C^1$ maps $u:N\lto M$. 
We define $\mN_k(N,M)$ as the sets of continuous  functions 
$F(t,x,v)$ such that
\begin{itemize}
\item $F/r_k$ is bounded.
\item $\int_NF(t,u(t),\d u_t)\d t=0$ for each $C^1$ map $u$.
\item There exists a compact   $K\subset N$ 
such that $F(t,x,v)=0$ if $t\not \in K$.
\item We have 
 $F(t,x,v)=\textsf{F}_{t,x}(v,\wedge_2 v,\cdots,\wedge_k v)=
 \textsf{F}(t,x,v,\wedge_2 v,\cdots,\wedge_k v)$, 
 where $\textsf{F}$ is continuous
 and 
where the functions $\textsf{F}_{t,x}( v,v_2,\cdots, v_k)$
are affine 
(We say that $F_{t,x}$ is $k$-affine).
\end{itemize}
By extension, we shall also denote by $\mN_k(N,M)$ 
the set of functions
$\textsf{F}(t,x,v_1,\cdots,v_k)$ associated with 
the elements $F\in \mN_k(N,M)$.
Note that the set $\mN_k(N,M)$ may depend on the metric on $M$ if $M$ is not compact.

\begin{defn}
A Young measure $\eta\in \mY_k$ is called
closed if $\int Fd\eta=0$ for all $F\in \mN_k$.
The set $\mC_k$ of closed measures is closed in 
$\mY_k$, and contains the Young measures $\hat u$
associated with maps $u\in W^{1,n}(N,M)$.
\end{defn}

Let us  explain how to build null-Lagrangians.
Given a field of $l$-vectors $U$, we will denote by $\dot U(t)$
the field of $(l-1)$-vectors which satisfies
$$
\d (i_U\Omega)=(-1)^{l+1}i_{\dot U} \Omega
$$
for any volume form $\Omega$ on $N$
which is compatible with the Riemannian metric
(meaning that the volume of an orthonormal base is $\pm 1$).
Notice that there are exactly two such volume forms 
on $N$ if it is orientable, and that they lead to the 
same $\dot U$. If $N$ is not orientable,
then no global volume form $\Omega$ exists,
but we can still define $\dot U$ by using volume forms defined
on orientable open subsets of $N$ (for example discs).
If $l=1$, for example, $U$ is a vector-field, and 
$\dot U=\text{div} U$.

\begin{lem}\label{null}
For each $l\in 1,\ldots, k$, each  smooth $(l-1)$-form $\chi$
on $M$ such that both $\chi$ and $\d\chi$ are bounded,
 and each compactly supported smooth field $U(t)$  of $l$-vectors on $N$, the function
$$
F(t,x,v):=
\chi_x\circ \wedge_{l-1}v\cdot \dot U(t)+
\d\chi_x\circ \wedge_l v\cdot U(t)
$$
belongs to $\mN_k(N,M)$.
In the case $l=1$, the form $\chi$ is just a function
$\chi(x)$ on $M$, and the function $F$ can be rewritten more clearly
$$
F(t,x,v)=\chi(x)\text{ \rm div} U(t)+\d \chi_x \circ v\cdot U(t).
$$
  \end{lem}

\proof
Let $u:N\lto \Rm^m$ be a $C^1$ function.
Let us still denote by $i_U\Omega$ and $\chi$
the pull-backs of $i_U\Omega$ and $\chi$
by the projections $N\times \Rm^m\lto N$ and
$N\times \Rm^m\lto \Rm^m$ respectively.
This allows us to define on $N\times \Rm^m$ the
$(n-1)$-form $\xi=i_U\Omega\wedge \chi$.
We have
$$
0=\int_N(Id\times u)^*\d \xi
=
\int_N\d\xi_{(t,u(t))}\circ\wedge_n(Id\times \d u_t)\cdot \lambda\d t
$$
$$
=(-1)^{l+1}\int_N (i_{\dot U}\Omega\wedge \chi)
\circ\wedge_n(Id\times \d u_t)\cdot \lambda \d t
+ (-1)^{n-l}
\int_N (i_{ U}\Omega\wedge \d\chi)
\circ\wedge_n(Id\times \d u_t)\cdot \lambda \d t
$$
Using Lemma \ref{calcul} of section 1.1,  we obtain
$$
0=
(-1)^{(l-1)(n-l)}
\int_N\chi_{u(t)}\circ \wedge_{l-1}du_t\cdot \dot U(t)\d t
+
(-1)^{(l+1)(n-l)}\int_N d\chi_{u(t)}\circ\wedge_l du_t \cdot U(t) \d t.
$$
After simplifying the signs, we obtain
$$
\int_N\chi_{u(t)}\circ \wedge_{l-1}\d u_t\cdot \dot U(t)\d t
+
\int_N d\chi_{u(t)}\circ\wedge_l \d u_t \cdot U(t) \d t=0.
$$
This is the required equality.
\qed

\subsection{Generalized maps and Cartesian maps}

The closed measure $\eta\in \mC_k$ is called a generalized map
if there exists a measurable map $u:N\lmto M$
such that the marginal of $\eta$ on $N\times M$
is concentrated on the graph of $u$.
We then say  that $\eta$ is a generalized map over $u$.
We denote by 
$\mG_k(N,M)$ the set of generalized maps.
 
\begin{defn}
We denote by $\ca(N,M)$ the set of measurable maps $u$ 
such that there exists a generalized map over $u$.
We call these maps cartesian maps.
We have a natural projection $\pi$ from
the set $\mG_k(N,M)$ of generalized maps to the set $\ca(N,M)$ of cartesian maps.
\end{defn}
The generalized maps have a remarkable structure:

\begin{thm}\label{structure}
Let $\eta$ be a generalized map over $u$.
Then, there exists a  measurable family $\Gamma_t$
of probability measures on $\mL(T_tN,T_{u(t)}M)$
such that
$\eta=dt\otimes \delta_{u(t)}\otimes \Gamma_t$.
Setting  
$$g_i(t):=\int_{\mL(T_tN,T_{u(t)}M)} \wedge_j v \;d\Gamma_t(v),
$$
 we have $u\in W^{1,1}(N,M)$, $g_1(t)= du_t$ and $g_i(t)=\wedge_i g_1(t)$
for almost all $t$.
\end{thm}
 By Jensen's inequality, we immediately obtain:

\begin{cor}\label{jensen}
If $\eta$ is a generalized map over $u$, and if $L$ is $k$-convex,
then
$$
\int_{J^1(N,M)}L \textrm{d}\eta\geq
\int_N L(t,u(t),du_t) dt. 
$$
\end{cor}
\proof
$$
\int_{J^1(N,M)}L \textrm{d}\eta=
\int_N \int_{\mL(T_tN,T_{u(t)}M)} L(t,u(t),v)\d \Gamma_t(v) \d t
$$
But we have, for each $t$,
$$
\int_{\mL(T_tN,T_{u(t)}M)} L(t,u(t),v)\d \Gamma_t(v)=
\int_{\mL(T_tN,T_{u(t)}M)} 
\textsf{L}(t,u(t),v,\wedge_2 v,\ldots,\wedge_k v)\d \Gamma_t(v)
$$
$$
\geq L(t,u(t),\d u(t),\wedge_2 \d u(t),\ldots,\wedge_k \d u(t))
$$
by Jensen's inequality, because 
$\int_{\mL(T_tN,T_{u(t)}M)}\wedge_j v\d \Gamma_t(v)=\wedge_j \d u(t)$ by Theorem \ref{structure}.
\qed

The proof of Theorem \ref{structure} will occupy the end of the present section.
The functions $g_i(t)$ depend only on the map $u$,
not on $\eta$. This is a consequence of the following:

\begin{lem}
 Let $u:N\lto M$ be a given measurable function.
 Then there exists at most one family of functions 
 $g_1(t), \ldots, g_k(t)$ 
 such that
 $$
 \int_N \textsf{F}(t,u(t), g_1(t), g_2(t),\ldots, g_k(t))dt =0
 $$
 for each $\textsf{F}\in \mN_k$.
  We call these functions 
 the distributional minors of $u$ if they exist. The map $u$
 belongs to $\ca(N,M)$ if and only if it admits distributional minors.
\end{lem}
\proof
The maps $g_l$ satisfy 
   the following equations:
 \begin{equation}\tag{E1}\label{E1}
\int_N \d\chi_{u(t)}\circ g_1(t) \cdot U(t) \d t+
\int_N \chi(u(t))\cdot \dot U(t)\d t=0 
\end{equation}
for all smooth vector-field $U$ on $N$ supported in the interior 
of $N$, and all smooth function $\chi:M\lto \Rm$,
and 
 \begin{equation}\tag{E\textit{l}}\label{El}
\int_N  \d\chi_{u(t)}\circ g_l(t) \cdot U(t) \d t+
\int_N \chi(u(t))\circ g_{l-1}(t)\cdot \dot U(t)\d t=0
\end{equation}
for all $l\in 2,\dots ,k$,  all
compactly supported  smooth field of $l$-vectors $U(t)$
on $N$, and  all smooth
$l-1$-form $\chi$ on $M$ which is bounded as well as $\d \chi$.
Now assume that $g'_l(t)$ are other maps satisfying the same equation.
Then, we have 
$$\int_N d\chi_{u(t)}\circ (g_l(t)-g'_l)\cdot U(t)=0$$
for each $l$, each $\chi$ and each $U$.
We claim that this implies that $g_l(t)-g'_l(t)=0$ almost
everywhere.
Since we have the freedom of choosing $U$, we conclude easily that
$d\chi_{u(t)}\circ (g_l(t)-g_l'(t))=0$
for almost all $t$.
If the claim 
 did not hold, we could find a compact set
$K\subset  N$ of positive measure, such that 
$u$ and $g_l-g_l'$  are continuous on $K$ and $g_l-g_l'$
does not vanish on $K$.
Let $t_0$ be a point of density of $K$, and let $\chi$
be a  compactly supported $(l-1)-$form on $M$ such that 
$$
d\chi_{u(t_0)}\circ(g_l(t_0)-g_l'(t_0))\neq 0.$$
Since $t_0$ is a density point of $K$, and since all the involved
functions are continuous on $K$,
there exists a compact subset $K'$ of $K$ of positive measure 
such that the relation
$
d\chi_{u(t)}\circ(g_l(t)-g_l'(t))\neq 0
$
holds for all $t\in K'$.
This is a contradiction.
\qed
\begin{lem}
If $u\in \ca(N,M)$, then $u\in W^{1,1}(N,M)$ and  the first distributional minor $g_1(t)$ of $u(t)$
is the weak derivative of $u$.
\end{lem}
\proof This is a direct consequence of (\ref{E1}).\qed

The following remark  can be applied for example when $f$ is an embedding
of $M$ into some Euclidean space, and $h$ is a chart of $N$:

\begin{prop}\label{charts}
Let $\tilde N$ and $\tilde M$ be other manifolds
and let $f:M\lto \tilde M$ and $h:\tilde N \lto N$.
Assume that $\tilde M$ is endowed with a complete metric.
If $f$ is smooth with bounded differential, $h$ is a smooth diffeomorphism onto its image
$h(\tilde N)\subset N$,
and $u\in \ca(N,M)$, then $f\circ u\circ h\in \ca(\tilde N,\tilde M)$.
Moreover, the distributional minors $\tilde g_i$ of $f\circ u\circ h$
are :
$$
\tilde g_i(\tilde t)= \wedge_i \d f_{u(h(\tilde t))}\circ g_i(h(t))\circ
\wedge_i \d h_t
$$
where $g_i$ are the distributional minors of $u$.
\end{prop}

\proof
Let us endow $\tilde N$ with the metric such that $h$ is an isometry.
Let $\tilde F(t,x,v)$ be an element of $\mN_k(\tilde N,\tilde M)$.
We want to prove that 
\begin{equation}\label{chart1}
\int_{\tilde N} \tilde F(\tilde t, f\circ u\circ h(\tilde t), 
\d f_{u(h(\tilde t))}\circ \d u_{h(\tilde t)}\circ dh_{\tilde t} )\d \tilde t =0.
\end{equation}
Setting 
$$
F(t,x,v)
:= \tilde F(h^{-1}(t),f(x),\d f_x\circ v\circ\d h_t)
$$
when $t\in h(\tilde N)$ and $F(t,x,v)=0$ when $t\not \in h(\tilde N)$,
we observe that (\ref{chart1}) is equivalent to
\begin{equation}\label{chart2}
\int_{ N} 
F(t,u(t),\d u_t)
\d t =0.
\end{equation}
This relation, on the other hand, holds if $F\in \mN_k(N,M)$
by definition of $\ca(N,M)$.

We prove that $F\in \mN_k(N,M)$. We begin to note that there is 
$\tilde K\subset\tilde N$, $\tilde K$ compact, such that 
$\tilde F(\tilde t,\tilde x,\tilde v)=0$ if $\tilde t\not\in\tilde K$; thus, 
$F(t,x,v)=0$ if $t$ does not belong to the compact set $h(\tilde K)$. Moreover, 
(\ref{chart2}) holds for all $C^1$ maps $u$.
This is true because (\ref{chart2}) is equivalent to
(\ref{chart1}), and  (\ref{chart1}) for $C^1$ maps follows because 
$\tilde F$ is a null Lagrangian, and thus it sends the $C^1$ map 
$f\circ u\circ h$ into zero; this amounts to (\ref{chart1}) by the chain rule.

In order to prove the equality between distributional minors,
we expand (\ref{chart1}) to 
$$
\int_{\tilde N} \tilde{ \textsf{F}}(\tilde t, f\circ u\circ h(\tilde t), 
\ldots,
\wedge_i \d f_{u(h(\tilde t))}\circ \tilde g_i(h(\tilde t))\circ \wedge_i \d h_{\tilde t},
\ldots 
)\d \tilde t =0
$$
and use lemma 7. 
\qed

\begin{lem}
We have $g_l(t)=\wedge_l g_1(t)$
for almost every $t\in N$.
\end{lem}
\proof
If $M$ is a Riemannian submanifold of the Euclidean space $E$,
then every map in $\ca(N,M)$ belongs to $\ca(N,E)$.
Therefore, using the embedding theorem of Nash, we can assume
for this proof that 
$M$ is a Euclidean space.
The set of points $t_0$ which are simultaneously Lebesgue points
of the function $u$ and of all the functions $g_l$,
have total measure.
Let $t_0$ be such a point.
By taking a chart in $N$, we can suppose that $N$ is the ball $B$ of radius
one in $\Rm^n$, that $t_0=0$, and that $\d t$ is  the Lebesgue measure. Translating in $\Rm^n$, we can suppose that $u(0)=0$.
Let us consider, for $s \geq 1$ the maps
$$
u^s(t):= su(t/s), 
\quad g_l^s:= g_l(t/s)
$$
on $N$.
By Proposition \ref{charts},  $u^s$ is a cartesian map on $\tilde N$, the ball of radius $s$, and 
$g_l^s$ are its distributional minors.
Our hypothesis on the point $t_0$
can be rephrased by saying that, strongly in $L^1(N)$, we have
$$u^s(t)\lto u^{\infty}(t)=g_1(0)t, 
\quad g^s_l(t)\lto g_l^{\infty}(t)= g_l(0)$$
when $s\lto \infty$. We can take a subsequence 
in order that these limits also hold almost everywhere.
Let $F$ be a null Lagrangian on the ball of radius $1$; in particular,
when trivially extended, it is a null Lagrangian on the ball of radius $s$, so that
$$\int_NF(t,u^s(t),\d u^s(t))\d t=0.$$
Passing to the limit, we obtain
$$
\int_N \textsf{F}(t,u^\infty(t),g_1^\infty(t),\cdots,g_k^\infty(t))\d t=0 .
$$
In other words, the limit function $u^\infty$ has $g_l^{\infty}$ 
as distributional minors.
On the other hand, since the function $u^{\infty}$
is smooth, we know that its distributional minors 
are $\wedge_ldu^{\infty}(t)$, which here are 
just the constant functions $\wedge _l g_1(0)$.
Therefore, by uniqueness of the distributional minors, 
we have  proved that 
$\wedge _lg_1(0)=g_l(0)$.
\qed

We have proved Theorem \ref{structure}. We can reformulate it
as follows: A function $u$ belongs to $\ca(N,M)$ if and only if 
the minors $\wedge_l du$ belong to $L^1$ and 
are distributional, which means that they satisfy the equation
$$
\int_N \textsf{F}(t,u(t),du_t,\wedge_2 du_t, \cdots,\wedge_k du_t)\d t=0.
$$
for all $\textsf{F}\in \mN_k$.
Note that $\text{cart}_1(N,M)=W^{1,1}(N,M)$.

\subsection{Topology}

The set $\mG_k$ of generalized maps is endowed with the topology of 
$\mY_k$.

\begin{prop}\label{closed}
The set $\mG_k(N,M)$ of generalized maps is closed in $\mY_k(N,M)$.
\end{prop}

\proof
Let $\eta^j$ be a sequence of generalized maps above $u^j$.
Let us assume that the sequence $\eta^j$ is converging to $\eta$
in $\mC_k(N,M)$. We have to prove that there exists a map 
$u\in W^{1,1}(N,M)$ such that the marginal of $\eta$ on $N\times M$
is concentrated on the graph of $u$.
It is enough to prove that, for each embedded ball $B\subset N$,
the marginal of $\eta_{|J^1(B,M)}$ on $B\times M$ is concentrated 
on the graph of a map $u$. As a consequence, we can suppose that 
$N$ is the open unit ball in $\Rm^n$.
We consider $M$ as a Riemannian submanifold of a Euclidean space $E$,
so that we see $u^j$ as elements of $W^{1,1}(N,E)$
with values in $M$.
Let $m^j\in E$ be the average of $u^j$, $m^j=\int_N u_j(t) \d t$.
Since the sequence $\eta^j$ is $r_k$-tight, see appendix,
the derivatives $\d u^j$ are bounded in $L^1$.
Therefore, by the Poincar\'e inequality, the sequence 
$(u^j-m^j)$ is bounded  in $W^{1,1}$.
By the compactness of the embedding $W^{1,1}\lto L^1$, 
this sequence is 
strongly compact in $L^1$. We assume, taking a 
subsequence,
that it has a limit $u^{\infty}$, and that the convergence holds
 almost everywhere.
By Lusin and Egorov Theorems, for all $\epsilon>0$, there exists a compact subset $K\in N$
such that $\d t(N-K)\leq \epsilon$ and such that 
$u^j$ is continuous on $K$ and
$(u^j-m^j)$ is  converging uniformly on $K$ to $u^{\infty}$.
It is clear at this point that the unboundedness of $m^j$
would contradict the tightness of $\eta^j$, and therefore
we can assume that the averages $m^j$ have a limit $m^{\infty}$.
Setting $u:=u^{\infty}+m^{\infty}$, we see that 
$u^j$ is converging uniformly to the continuous function $u$
on $K$.
Denoting by $\mu$ the marginal of $\eta$ on $N\times M$,
we conclude that the $\mu$-measure of the graph of $u$
is greater that $1-\epsilon$. Since this holds for all $\epsilon>0$,
we conclude that the measure $\mu$ is concentrated on the graph 
of $u$.
\qed

\section{Boundary conditions and compactness}\label{compactsect}
In most applications, the manifold $M$ is not compact, 
and it is necessary to introduce boundary conditions in order
to get compactness.
We fix, as explained in the introduction,
a map $u_0\in \ca(N,M)$.
We define the set $ \mC_k(N,M;u_0)\subset \mC_k(N,M)$
 of closed measures
with boundary $u_0$ as
the set of
measures $\eta\in \mC_k(N,M)$ such that 
$$
\int_{\J}\d \chi_x\circ v\cdot U(t)+\chi(x)\cdot 
\text{\rm div}U(t)\d \eta  (t,x,v)
=
\int_N \d \chi_{u_0(t)}\circ \d u_0(t)\cdot U(t) +
 \chi(u_0(t))\cdot \text{\rm div}U(t)\d t.
$$
for each smooth vectorfield $U(t)$ on $\bar N$ (not necessarily  supported in a compact set of $N$)
and each bounded smooth function $\chi(x)$ on $M$ with bounded derivative.
We can also 
define the set of generalized maps with boundary value $u_0$:
$$
\mG_k(N,M;u_0):=\mG_k(N,M)\cap \mC_k(N,M;u_0).
$$
The space $\ca(N,M;u_0)$ 
is  the space of maps $u$ such that there exists
a generalized map
$\eta\in \mG_k(N,M;u_0)$ above $u$, or in other words the maps 
$u$ whose Young measure  $\hat u$ belongs to $\mG_k(N,M;u_0)$. 
In \cite{GMS1}, the functions in $\ca(N,M;u_0)$ are said to satisfy a weak anchorage condition.
In the case where $M$ is a Euclidean space, the functions
in $\ca(N,M;u_0)$ are just the functions in $\ca(N,M)$
which have the same trace on $\partial N$ as $u_0$ in the 
$W^{1,1} $ sense. 

\begin{prop}\label{compact} 
Let $L$ be a $k$-convex and coercive Lagrangian.
For each $c>0$, the set of measures $\eta\in \mC_k(N,M;u_0)$
which satisfy 
\begin{equation}\label{bound}
\int_{J^1(N,M)}L(t,x,v){\rm d}\eta(t,x,v)\leq c  
\end{equation}
is compact.
\end{prop}

\proof 
Let us denote by $\mC(c)$ the set of measures $\eta\in \mC_k(N,M;u_0)$
which satisfy (\ref{bound}).
Since the functional
$$\eta\lmto
\int Ld\eta
$$
is lower semi-continuous on $\mY_k(N,M)$ (by Proposition \ref{lsc}), 
and since $\mC_k(N,M;u_0)$, is closed in $\mY_k(N,M)$,
the set $\mC(c)$ is closed in $\mY_k(N,M)$.
So it is enough to prove that it is relatively compact.
By the appendix, this follows if we can prove that it
is $r_k$-tight.
In other words, we have to show 
that  for each $\epsilon>0$ there exists a compact subset $Z\in \J$
such that 
$$
\int_{\J \setminus Z(R)}r_k(t,x,v){\rm d}\eta(t,x,v)\leq 2\epsilon
$$
for each measure $\eta\in \mC(c)$.
We shall prove that this holds for 
$$Z(R)=\{ (t,x,v)\in \J \quad\colon\quad d(x_0,x)\leq R,||v||\leq R \}  
 $$
 when $R$ is large enough
($x_0$ is a point in $M$ that we have fixed once and for all).
At this point it is convenient to assume, without loss of generality,
that $L\geq 0$.
We define 
$$\tilde Z(R)=\{ (t,x,v)\in \J \quad\colon\quad ||v||\leq R \} $$
and we see that there exists $A(R)>0$ with $A(R)\rightarrow+\infty$ as $R\rightarrow+\infty$
such that, for all $\eta\in \mC(c)$, 
$$
c\ge\int_{\J\setminus\tilde Z(R)}
L(t,x,v){\rm d}\eta(t,x,v)\geq A(R)
\int_{\J\setminus\tilde Z(R)}r_k(t,x,v){\rm d}\eta(t,x,v).
 $$
Taking $R$ 
sufficiently large, we get from the inequality above that
\begin{equation}\label{speed}
\int_{\J\setminus\tilde Z(R)}
r_k(t,x,v)\d\eta(t,x,v)\leq{\epsilon}
\end{equation}
for each $\eta \in \mC(c)$.
Setting now 
$$\hat Z(R)=\{ (t,x,v)\in \J\quad\colon\quad d(x_0,x)\leq R \} $$
we see that  the desired inequality  follows if we prove that
\begin{equation}\label{space}
\eta(\J-\hat Z(R))\leq\epsilon(R)
\qquad \forall \eta\in \mC(c)
\end{equation}
for all  $R$, with $\epsilon (R)\lto 0$ as $R\lto \infty$.
Indeed, taking $R_0$ such that (\ref{speed}) holds, and then
setting $S=\max _{\tilde Z(R_0)}r_k$, we get
$$
\int_{\J\setminus Z(R)}r_k \d \eta\leq 
\int_{\J \setminus \tilde Z(R_0)}r_k \d \eta
+S\eta(\J\setminus \hat Z(R))\leq 
\epsilon + S\epsilon(R).
$$
In order to prove (\ref{space}), we consider, for each $R>0$, a 
function 
$g\in C^1(M,\Rm)$ such that
$$0\leq g_R(x)\leq 1,\quad g(x)=1\quad\hbox{if}\quad d(x_0,x)
\geq R,$$
$$
g(x)=0\quad\hbox{if}\quad d(x_0,x)\leq R/2
$$
$$
|d g_x|\leq \delta(R)\quad\forall x  
  $$
  where $\delta(R)\lto 0$ as $R\lto \infty$;
  and a smooth vector-field $U(t)$ on $\bar N$  such that 
  $\dot U=1$ on $N$
  or equivalently such that $\text{div } U=1$ on $N$.
  The existence of such a vector-field is given by Lemma 
  \ref{conform} below. We note that $U$ is bounded, since $\bar N$ is compact.
 Since $\mC(c)\subset \mC_k(N,M;u_0)$, we  have, for $\eta\in \mC(c)$,
$$
\int_{\J} g(x) {\rm d}\eta(t,x,v)=
\int_{N} g(u_0(t))+\d g_{u_0(t)}\circ \d u_0(t)\cdot U(t)\d t
-\int_{\J}{\rm d}g_x\circ v\cdot U(t) {\rm d}\eta(t,x,v).
$$
The last formula and the 
definition of $g$ imply that there exists $C>0$ such that
$$
\int_{\J \setminus Z(R)}{\rm d}\eta(t,x,v)\leq
C\delta(R)+\int_{N} g(u_0(t))\text{d}t
$$
for all $R$ and all $\eta \in \mC(c)$.
The term on the right converges to zero as $R\lto \infty$,
this ends the proof.
\qed

\begin{lem}\label{conform}
Let $\bar N$ be a compact Riemannian manifold with a non-empty boundary. 
There exists a smooth  vector-field $U(t)$ on 
$\bar N$ such that $\dot U= 1$ on $N$ 
 or equivalently such that 
$\text{\emph{div} } U= 1$ on $N$.
\end{lem}
\proof
In the case where $ N$ is a ball in $\Rm^n$, this is obvious,
just  take $U(t)=t/n$.
In general, one can build $U$ as the gradient of a function
$h$ which solves $\Delta h=1$ on $N$.
\qed

\section{Conclusion}\label{conclusion}
We now collect the tools we have introduced to prove Theorem
\ref{main}. We also add some discussions and variations.

\subsection{Proof of Theorem \ref{main}}

By Propositions \ref{closed} and \ref{compact},
there exists a generalized map $\eta$ over some $u\in W^1(N,M)$ such that $\eta$ minimizes $\int L\d \eta$ on $\mG_k(N,M;u_0)$. We want to show that $u$ minimizes in $\ca(N,M;u_0)$.

If $v\in \ca(N,M;u_0)$ is another map,
we have
$$
\int_{J^1(N,M)} L(t,u(t),\d u_t)\d t
\leq \int L \d\eta
\leq \int L \d \hat v
=\int_{J^1(N,M)} L(t,v(t),\d v_t)\d t
$$
where the first inequality comes from Corollary \ref{jensen},
and where $\hat v\in \mG_k(N,M;u_0)$ is the Young measure associated with $v$.
This proves that $u$ is minimizing in $\ca(N,M;u_0)$.

\subsection{An example}\label{example}
We consider $M=\Rm^2$; $N=B$,  the unit
open ball of 
$\Rm^2$, and the Lagrangian 
$$
L(t,x,v)=\epsilon (|v|^p+|t|^4|v|^4 )+|\text{det } v|^2,
$$
with $p\in ]1,2[$.
We claim that 
$$
\inf _{u\in W^{1,1}(N,M;Id)} \int_N L(t,u(t),\d u_t)\d t
<
\min _{u\in \text{\rm cart}_2(N,M;Id)}\int_N L(t,u(t), \d u_t)dt
$$
when $\epsilon>0$ is small enough.
Indeed, taking $\underline u(t)=t/|t|$, and observing that
 $\det \d \underline u =0$, we get a constant $C>0$ such that
$$
\inf _{u\in W^{1,1}(N,M;Id)} \int_N L(t,u(t),\d u_t)\d t
\leq \int_N L(t,\underline u(t),\d \underline u_t)\d t
\leq C\epsilon.
$$
On the other hand, if $u\in \text{\rm cart}_2(B,\Rm^2;Id)$ is a minimizer we have 
$$
\int_B L(t,u(t),\d u_t)\d t\geq
\int_B|\det\d u_t|^2\d t\geq
{1\over |B|}\big(\int_B|\det\d u_t|\d t\big)^2
\geq \pi,
$$
where the last inequality follows from the following Lemma:

\begin{lem}
If $u\in \text{\rm cart}_2(B,\Rm^2;Id)$ minimizes
$\int_B L (t,u(t),\d u_t)\d t<\infty$,
then 
$\int_B \det \d u (t) \d t=\pi$.
\end{lem}
\proof
We claim that $u(B)\subset \bar B$.
Indeed, let $f:\Rm^2\lto \Rm^2$ be a smooth diffeomorphism such
that $f=Id$ on $B$ and $|\d f|<1$ outside of $\bar B$.
By Proposition \ref{charts}, the map $f\circ u$ belongs to $\text{cart}_2(B,\Rm^2)$,
and it has the same boundary condition as $u$.
Since $|\d f(u)|\leq 1$ we have that $|\d (f\circ u)|\leq|\d u|$ and 
$|\det\d (f\circ u)|\leq |\det\d u|$; the first inequality is strict if $|u|>1$ and $\d u\not=0$. 
If we did not have $u(t)\subset \bar B$ for almost every $t$, the action of $f\circ u$
would be strictly smaller than the action of $u$, which would contradict
the assumption that $u$ is a minimizer.

Let us denote by $A$ the annulus $1/2<|t|<1$.
We have 
$$u\in W^{1,p}(B,\Rm^2)\cap W^{1,4}(A,\Rm^2),
$$
so that $u$ is continuous on $A$, and extends by continuity
to $\partial B$, where it takes the value $u|_{\partial B}=Id$.
Finally, recall that $u(B)\subset \bar B$.
Define
$$
u_i(t):={i^2}\int_B \tau(is)u((1-1/i)t-s)\d s 
$$
where $\tau:B\lto [0,1]$ is a smooth convolution kernel.
It is classical that $u_i\lto u$ in $W^{1,p}(B,\Rm^2)$,
and in $W^{1,4}(A,\Rm)$. As a consequence, 
$u_i|_{\partial B}$ converges uniformly to the identity; since $u_i$ is smooth, this implies
$$\int_B \det \d u_i \d t\lto \pi.
$$
Thus it is enough to prove that
\begin{equation}\label{spirit}
\int_B\det\d u_i\d t\lto \int_B\det\d u\d t . 
\end{equation}
Let $r\in]1/2,1[$ and let $\phi_r\in C^\infty_0(B,\Rm)$ be such that 
$0\leq\phi_r\leq 1$ and $\phi_r=1$ on $B(0,r)$. Since 
$$\int_B\det\d u\cdot\phi_r\d t\lto \int_B\det\d u\d t$$
as $r\lto 1$, the formula (\ref{spirit}) follows if we prove
\begin{equation}\label{flic}
\int_B\det\d u_i\cdot\phi_r\d t\lto \int_B\det\d u\cdot\phi_r\d t 
\quad\forall r\in]1/2,1[ 
\end{equation}
and
\begin{equation}\label{floc}
\left\vert \int_B\det\d u_i\cdot(1-\phi_r)\d t \right\vert <\epsilon
\quad\forall i
\quad\text{if}\quad r\ge 1-\delta . 
\end{equation}
Note that (\ref{floc}) follows  from the boundedness of $\det \d u_i$ in $L^2(A)$.
In order to prove (\ref{flic}),
we set $u_i=(u^1_i,u^2_i)$ and $u=(u^1,u^2)$, we call $(x_1,x_2)$ 
the coordinates on the target space $\Rm^2$ and we 
assert that
$$
\int_B\det\d u\cdot\phi_r\d t=
\int_B  
u^1(\partial_2u^2,-\partial_1u^2)    \cdot
\nabla\phi_r\d t  . 
$$
Indeed, this formula is just (\ref{El}) with $l=2$, 
$\chi=a(x^1)\d x^2$ and $U=\phi_r e_1\wedge e_2$,
where $a:\Rm\lto [-3,3]$ is a smooth function such that $a(x^1)=x^1$ on $[-2,2]$.
Here we use that $u(B)\subset \bar B$.
Similarly, by Lemma \ref{null},
$$\int_B\det\d u_i\cdot\phi_r\d t=
\int_B  
u^1_i(\partial_2u^2_i,-\partial_1u^2_i)    \cdot
\nabla\phi_r\d t  .  $$
As a consequence (\ref{flic}) is equivalent to 
$$
\int_B
[u^1_i(\partial_2u^2_i,-\partial_1u^2_i)]\cdot
\nabla\phi_r\d t\lto 
\int_B
[u^1(\partial_2u^2,-\partial_1u^2)]\cdot
\nabla\phi_r\d t
   $$
which holds because the integrand is converging almost everywhere
 and is 
bounded in $L^2$.
 \qed

\subsection{Weak continuity of minors}
Let us mention the following classical result
which follows from our tools (see \cite{GMS1}, 3.3.1
or \cite{Da}, 8.3):

\begin{prop}
Let $N$ be a bounded disc in $\Rm^n$.
Let $u_i$ be a sequence of maps in $\ca(N,\Rm^m)$,
and let $u\in W^{1,1}(N,M)$ and 
$g_j(t)\in L^1(N,\mL(\wedge_j\Rm^n,\wedge_j\Rm^m))$
be   such that $u_i\lto u$ weakly in $W^{1,1}$ and 
$$ \wedge _2 \d u_i\lto g_2
,\quad\ldots\quad \wedge _k \d u_i\lto g_k
$$
weakly in $L^1$.
Then $g_2=\wedge _2\d u,\ldots g_k=\wedge _k\d u$.
\end{prop}
\proof
We consider the Young measures 
 $\hat u_i$ in $\mG_k(N,M)$ associated with the functions $u_i$.
Now weak convergence implies uniform integrability, which translates to the fact that $\hat u_i$ is $r_k$-tight,
 and therefore compact in $\mY_k(N,\Rm^m)$.
 We can suppose that it has  a limit $\eta$,
 which is a generalized map above $u$.
 If 
 $F(t,v)=\textsf{F}(t,v,v_2,\ldots,v_k)$ 
 is a 
 continuous function which is affine in $(v,v_2,\ldots,v_k)$,
 then we have 
 $$\int F\d \hat u_i\lto \int F\d \eta$$
 because $\hat u_i\lto \eta$.
 On the other hand, since $F$ is affine in the minors, 
 and since $\wedge _j \d u_i\lto g_j$ weakly, we have
 $$
 \int \textsf{F}(t,\d u_i(t),\wedge_2 \d u_i(t)
 ,\ldots,\wedge_k \d u_i(t))\d t\lto 
 \int \textsf{F}(t,\d u(t), g_2(t),\ldots g_k(t))dt.
 $$
  We conclude that 
 $$
 \int F\d \eta =\int \textsf{F}(t,\d u(t), g_2(t),\ldots g_k(t))dt.
$$
This implies that, for almost all $t$,
$$g_j(t)=\int \wedge_j v \d\Gamma_t=\wedge_j \d u(t)$$
by Theorem \ref{structure}.
\qed

\subsection{On Null-Lagrangians}
It may seem unnatural in the definitions of the sets 
$\mN_k(N,M)$ to require that the null-Lagrangians
$F(t,x,v)$
be $k$-affine functions of $v$.
Indeed, working with a larger set $\mN(N,M)$
of null-Lagrangians would make the result stronger,
and may allow to relax somewhat the $k$-convexity hypothesis on $L$.
The following result, however, shows that there is not much hope in that
direction:

\begin{prop}
Let us assume that $N= \Rm^n$, and that $M=\Rm^m$.
If $F$ is a null-Lagrangian  such that $F/r_k$ is bounded,
then $F_{t,x}$ is $k$-affine for each $t$ and $x$.
 \end{prop}
\proof
We just give an idea of the proof.
It follows from proposition 9 that, for $\lambda >0$
and $(t_0,x_0)\in \Rm^n\times \Rm^m$,
the function
$$
F_{\lambda}(t,x,v):= F(t_0+\lambda t, x_0+\lambda x, v)
$$
is  a null-Lagrangian.
But then $F_0$ is also a null-Lagrangian,
which means that
 $F_{t_0,x_0}$ is quasi-affine in the sense of \cite{Da}, section 4.1; but in the same section of \cite{Da} it is proven that  quasi-affine
functions  are  poly-affine.
In other words, 
there exists an affine  function $\mathsf{F}(v,v_2,\ldots, v_p)$
such that $F_{t_0,x_0}(v)=\mathsf{F}(v,\wedge_2 v,\ldots,\wedge_p v)$,
where $p=\min\{m,n\}$.
But the bound implies that $\mathsf{F}$ does not depend on 
$\wedge_j v$ for $j> k$. 
\qed

\subsection{More general setting}
The heart of the matters is the Jensen's inequality obtained
in Corollary \ref{jensen}.
This inequality is the result of an equilibrium between
the known properties of the measures $\Gamma_t$
appearing in the disintegration of generalized maps
and the convexity assumed on the integrand $L$.

Other, but less explicit equilibria might be obtained as follows.
Let $r(t,x,v)$ be a continuous function
on $\J$ such that $r(t,x,v)\geq 1+ \|v\|$.
We define the associated Kantorovich-Rubinstein space
$\mP_r(\J)$, which is the set of Borel probabilities $\eta$
on $\J$ such that $\int r\d \eta< \infty$.
We also define the set $\mY_r(N,M)$ of those elements $\eta$
of $\mP_r(\J))$
such that $t_{\sharp} \eta=\d t$.

Now let $\hat \mG_r(N,M)$ be the closure, in  
$\mY_r(N,M)$ of the set of Young measures associated with 
smooth maps.
We can prove as in Proposition \ref{closed} that, to each 
$\eta\in \hat \mG_r(N,M)$ is associated a map
$u\in W^{1,1}(N,M)$ such that 
\begin{equation}\label{dec}
\eta=\d t\otimes \delta_{u(t)}\otimes \Gamma_t
\end{equation}
and such that $\int v\d \Gamma_t=\d u(t)$ for almost all $t$.
We can define $\text{cart}_r(N,M)$ as the set of maps which appear 
in this way.
In this setting, we can fix boundary conditions
as before by taking $u_0\in \text{cart}_r( N,M)$.
If the coercivity condition of the Lagrangian is modified to
$$
L(t,x,v)\geq l(r(t,x,v)),
$$
with $l$  super-linear, we still have compactness:
Proposition \ref{compact} still holds, with the same proof.
So if $L$ is a normal integrand satisfying the modified coercivity condition,
then there exists a Young measure $\eta\in\hat \mG_r(N,M;u_0)$
which minimizes the integral $\int Ld\eta$ in this set.

In order to prove the existence of minimizers in 
$\text{cart}_r(\tilde N,M;u_0)$, it is enough to adapt the 
convexity condition, in such a way that Corollary 
\ref{jensen} holds for the elements of $\hat \mG_r(N,M)$.

Let $\mP_{t,x}$
be the set of Borel probability measures $\Gamma$ on $\mL(T_tN,T_xM)$
such that $\int r_{t,x}(v)d\Gamma(v)<\infty$.
In short, we have
 $$\mP_{t,x}:=\mP_{r_{t,x}}(\mL(T_tN,T_xM))
 $$
(see the Appendix below).
Let $B$ be a closed  ball of volume one in $T_tN$.
Let $\mathsf{P}_{t,x}$ be the closure, in $\mP_{t,x}$
of the measures of the form
$$\Gamma=(\d u)_{\sharp} (\d t_{|B})
$$
where $u:T_tN\lto T_xM$ is a smooth map supported in $B$.
Note that if $\Gamma\in \mathsf{P}_{t,x}$, then 
$\int_{\mL(T_tN,T_xM)} v \d  \Gamma =0$. 
A last notation is necessary: we denote by $\tau_z$ the translation
of vector $z$.
Then, possibly under some mild assumption on the function r,
the following result can be proved by a blow-up argument called localisation procedure in \cite{Pe}:
\vs

\noindent
\textbf{Structure Theorem:}
\begin{itshape}
The measures $\eta\in \hat \mG_r(N,M)$
can be written in the form (\ref{dec}), with 
$$
\left (\tau_{-\d u_t }\right )_{\sharp}\Gamma_t\in \mathsf{P}_{t,x}$$
for almost all $t$.
\end{itshape}\vs

As a consequence, the convexity condition that has to be assumed in order
that Corollary \ref{jensen}, and then Theorem \ref{main}
hold in this more general setting is 
$$
\int_{\mL(T_tN,T_xM)}L_{t,x}(a+v)\d\Gamma(v) \geq L(t,x,a)
$$
for all $(t,x)\in N\times M$, for all $a\in  \mL(T_tN,T_xM)$
and for all $\Gamma\in \mathsf{P}_{t,x}$.
This is not an easy condition to check on examples.

\begin{appendix}
\section{Kantorovich-Rubinstein space}
Let us recall some standard facts on probability
measures, see \cite{AGS,Vi:03}. 
 Let $(X,d)$ be a complete and separable metric space,
 and let $r:X\lto [1,\infty)$ be a continuous function.
 Let $\mP_r(X)$ be the set of Borel probability measures 
 $\mu$ on $X$ 
 which satisfy
 $$\int_X r(x)  d\mu(x)<\infty.$$
 Let us denote by $C_r(X)$ the set of continuous functions $f$  on $X$
such that 
$$
\sup_{x\in X} \frac{|f(x)|}{r(x)}<\infty.
$$
There exists a distance $d$ on $\mP_r(X)$ such that 
$d(\mu_n,\mu)\lto 0$
if and only if
$$
\int fd\mu_n\lto \int fd\mu
$$
for all $f\in C_r(X)$.
This distance can be chosen such that, in addition,
the metric space $(\mP_r,d)$ is a complete and separable metric space.

In order to define such a distance $d$ on $\mP_r(X)$
one can define first the distance 
 $$
d_r(x,y):=\min(d(x,y),1)+|r(y)-r(x)|.
$$
on $X$, which is complete and equivalent to $d$.
Then, we can define the distance $d$ on $\mP_r(X)$
as the Kantorovich-Rubinstein
(also called 1-Wasserstein) distance of $(X,d_r)$.

The relatively compact subsets of $(\mP_r(X),d)$ are those which are 
$r$-tight:

 \begin{defn}
 The subset $Y\subset \mP_r(X)$
 is called $r$-tight if one of the following equivalent properties holds:
 \begin{itemize}
 \item For each $\epsilon>0$, there exists a compact set
 $K\subset X$  such that
 $\int_{X-K} r(x)d\mu \leq \epsilon$
 for each $\mu \in Y$.
 \item There exists a  function $f:X\lto [0,\infty]$
 whose sublevels are compact
 and  a constant $C$ such that 
 $\int_X r(x)f(x)d\mu\leq C$ for each $\mu \in Y$.
\item The family $Y$ is tight and $r$ is $Y$-uniformly integrable.
The first means that, for each $\epsilon>0$, there exists a compact set
 $K\subset X$  such that
 $\mu (X-K) \leq \epsilon$
 for each $\mu \in Y$.
 The second means that for each $\epsilon>0$, there exists
 a ball $B$ in $X$ such that 
 $\int_{X-B} r(x)d\mu \leq \epsilon$
 for each $\mu \in Y$.
 \end{itemize}
\item
\end{defn}
 Note that $1$-tightness is just tightness if $r\equiv 1$.
 If $r$ is proper, then $Y$ is $r$-tight if and only if 
there exists a constant $C$ and a superlinear function
$f:[0,\infty)\lto \Rm$ such that
$$
\int _X f\circ r d\eta\leq C
$$
for all $\eta\in Y$.

\end{appendix}

\small

\end{document}